# Artificial Neural Network Approach for Solving Fractional order initial value problems

Susmita Mall and S. Chakraverty*

Department of Mathematics,

National Institute of Technology Rourkela, Odisha-769008, India

Tel: 91661-2462713, Fax: +91661-2462713-2701

E-mail*: sne_chak@yahoo.com

## Abstract

In this paper, an Artificial Neural Network (ANN) technique is developed to find solution of celebrated Fractional order Differential Equations (FDE). Compared to integer order differential equation, FDE has the advantage that it can better describe sometimes various real world application problems of physical systems. Here we have employed multi-layer feed forward neural architecture and error back propagation algorithm with unsupervised learning for minimizing the error function and modification of the parameters (weights and biases). Combining the initial conditions with the ANN output gives us a suitable approximate solution of FDE. To prove the applicability of the concept, some illustrative examples are provided to demonstrate the precision and effectiveness of this method. Comparison of the present results with other available results by traditional methods shows a close match which establishes its correctness and accuracy of this method.

# Introduction

Recently, fractional calculus has become popular in the scientific community because it has numerous applications in various fields of science and engineering [1-5]. Fractional differential equations describe various phenomena such as fluid flow in porous material, anomalous diffusion transport, signal processing, control theory of dynamical systems, viscoelasticity etc. [7]. In most cases, analytical solutions of fractional order differential equations may not be obtained easily. Therefore, it is important to develop some reliable and efficient techniques to handle FDEs.

During the past decades, various numerical techniques have been developed to solve these equations. These methods include finite difference [8-10], Adomian decomposition [11, 12], Wavelet method [13], variational iteration [14, 15], Laplace transforms [16, 17] and operational matrix method [18-22] etc. Although these techniques provide good approximations to the solution, but most of these methods give series expansions in the neighborhood of initial conditions are used [23]. Besides that, rounding-off errors solemnly influence the solution precision in the numerical methods with complicatedness that also increase quickly with the number of sampling points [6].

In last few years, various machine intelligence procedures in particular Artificial Neural Network (ANN) methods have been established as a powerful technique to solve a variety of real-world problems because of its excellent learning capacity [24, 25]. ANN approach has attracted much consideration to its advantages such as learning, adaptive, error computation, fault-tolerance etc. Currently, a lot of attention has been devoted to the study of ANN for solving ordinary and partial differential equations [26-31]. Here our target is to solve fractional order initial value problems using ANN.

According to the historical survey, artificial neural network based approaches are utilized by various researchers to solve different types of fractional order differential equations. In [6] Pakdaman et al. employed neural network and Broyden–Fletcher–Goldfarb–Shanno (BFGS) optimization technique to solve linear and nonlinear FDEs. Jafarian et al. [32] have applied artificial neural network model for approximate polynomial solution of special type of fractional order Volterra integro differential equations. The neural network training and cosine basis functions with adjustable parameters have been presented by Qu and Liu [33] for solving single and the systems of coupled fractional order differential equations. Rostami and Jafarian [34] used the combination of ANN approach and power series method for handling higher order linear fractional differential equations. In [35], Almarash implemented an approximation solution of fractional partial differential equations by using neural network method with radial basis functions. Sabouri et al. [36] employed the ability of neural networks for solving fractional order optimal control problems. In another study [37] Jafarian et al. have used combination of multi-layer ANN and the power series method to the numerical solution of a class of fractional order initial value problems.

As per the review of the literatures it reveals that the authors have used multi-layer ANN with optimization technique for solving FDEs [6]. Previously, few researchers have considered three layer feed forward ANN structure and power series method for the numerical solution of different types of fractional order initial value problems [32, 35 and 37]. They have taken the solution function as a series polynomial in which its coefficients are determined by power series method.

In this investigation, the authors vested their effort to develop a multi-layer ANN model with unsupervised back propagation learning algorithm for solving fractional order initial value problems. The ANN approximate solution of FDE is written as sum of two terms, first one satisfies initial/boundary conditions and the second part involves output of neural network with adjustable parameters. Feed forward neural network model and error back propagation principles (gradient descent procedure) are used for modification of the network parameters and to minimize the computed error function. Initial weights from input to hidden and hidden to output layer are considered as random. The approximate solution of FDEs by ANN is found to be advantageous but it depends upon the ANN model that one considers such as

- Solution search proceeds without coordinate transformations;
- Simple implementation and easy computation;
- The back propagation algorithm is unsupervised;
- After training of the ANN model, we may use it as a black box to get numerical results at any arbitrary point in the domain etc.

Some examples are also given to show the efficiency of the procedure.

We begin our work in section 2, by reviewing certain basic definitions and fundamental issues of fractional calculus. ANN formulations for FDEs, construction of the appropriate form of the approximate solution of fractional order initial value problem and error estimation are described in section 3. In section 4, we have presented the numerical examples, solutions and comparison of analytical and ANN results. Lastly section 5, incorporates the conclusions.

## 2. Preliminaries

In this section, we recall some definitions and general concepts related to fractional calculus which may be used further in this paper [38-41]

***Definition 2.1: Riemann-Liouville type fractional derivative* [39]**

The Riemann-Liouville type fractional derivative of order $\alpha > 0$ of a function $f:(0,\infty) \to R$ is given as

$$D^{\alpha} f(t) = \frac{d^n}{dt^n} \frac{1}{\Gamma(n-\alpha)} \int_0^t (t-\tau)^{n-\alpha-1} f(\tau) d\tau \tag{1}$$

Where $n = |\alpha| + 1$ and $|\alpha|$ is the integer part of $\alpha$.

***Definition 2.2: Riemann-Liouville type fractional integral* [38]**

The Riemann-Liouville type fractional integral of order $\alpha > 0$ of a function $f:(0,\infty) \to R$ is

$$I^{\alpha} f(t) = \frac{1}{\Gamma(\alpha)} \int_0^t (t-\tau)^{\alpha-1} f(\tau) d\tau \tag{2}$$

Here $\Gamma$ denotes the Gamma function.

## Definition 2.3: Caputo type fractional derivative [39, 40]

The Caputo type derivative of order $\alpha$ and $\alpha \in [n-1, n)$, is expressed as

$$D_a^\alpha f(t) = \frac{1}{\Gamma(n-\alpha)} \int_a^t \frac{f^n(\tau)}{(t-\tau)^{\alpha-n+1}} d\tau \tag{3}$$

## Definition 2.4: Conformable fractional derivative [42, 43]

Let us consider a function $f : [0, \infty) \to R$ and $\alpha \in (0,1]$. Then the conformable fractional derivative of $f$ of order $\alpha$ is defined as

$$T_\alpha(f(t)) = \lim_{\varepsilon \to 0} \frac{f(t + \varepsilon t^{1-\alpha}) - f(t)}{\varepsilon} \tag{4}$$

For all $t > 0$ and $T_\alpha(f(t))$ denotes the conformable fractional derivative of order $\alpha$.

$T_\alpha(f(t))$ satisfies all the following properties

I. $T_\alpha(t^p) = pt^{p-\alpha}$     for all $p \in R$

II. $T_\alpha(\eta) = 0$     for all constant function $\eta$

III. $T_\alpha(fg) = fT_\alpha(g) + gT_\alpha(f)$

IV. $f$ is differentiable, then $T_\alpha(f(t)) = t^{1-\alpha} \frac{df}{dt}(t)$

V. $T_\alpha\left(\frac{1}{\alpha} t^\alpha\right) = 1$

VI. $T_\alpha(e^{cx}) = cx^{1-\alpha} e^{cx}$,     $c \in R$

VII. $T_\alpha(\sin bx) = bx^{1-\alpha} \cos bx$     etc.

In this paper, we have used the properties of Conformable fractional derivative.

## 3. Description of the Method

In this section, we describe general formulation of ANN model for fractional order differential equation. In particular, structure of ANN model and the formulations for initial value fractional order problems are incorporated in detail.

### 3.1 ANN formulation for fractional order differential equations

Let us consider a fractional order initial value problem as

$$(_{x_0}D_x^\alpha y(x)) = f(x, y(x)) \tag{5}$$

subject to $y(x_0) = y_0$

where $n-1 < \alpha \leq n$, $n \in N$

Let $y_N(x, \Omega)$ denotes the approximate solution of ANN model with $\Omega$ is a vector containing corresponding weights and $x$ is the input data. The above FDE is transformed into the following problem

$$(_{x_0}D_x^\alpha y_N(x, \Omega)) = f(x, y_N(x, \Omega)) \tag{6}$$

The approximate solution $y_N(x, \Omega)$ of feed forward neural network with network parameters $\Omega$ may be written in the form

$$y_N(x, \Omega) = A(x) + F(x, N(x, \Omega)) \tag{7}$$

The first term $A(x)$ in right hand side does not contain adjustable parameters and satisfies only initial\boundary conditions, whereas the second term $F(x, N(x, \Omega))$ contain the single output $N(x,\Omega)$ of feed forward neural network with input $x$ and vector containing the corresponding weights $\Omega$.

Here we consider a three layer ANN model with one input node $x$, one hidden layer consisting of $m$ number of nodes and one output node $N(x,\Omega)$.

The output $N(x,\Omega)$ is expressed as

$$N(x,\Omega) = \sum_{j=1}^{m} v_j \varphi(z_j) \tag{8}$$

where $z_j = w_j x + u_j$ and $w_j$ is weight from input to $j^{th}$ hidden unit, $v_j$ denotes the weight from $j^{th}$ hidden unit to output unit and lastly $u_j$ is the bias for $j^{th}$ hidden node. In this investigation we have considered the sigmoid function $\varphi(x) = \dfrac{1}{1+e^{-x}}$ as an activation function.

General form of corresponding error function for the fractional order initial value problem may be formulated as

$$E(x,\Omega) = \sum_{i=1}^{h} \left\{ (_{x_0} D_x^{\alpha} y_N(x_i,\Omega)) - f(x_i, y_N(x_i,\Omega)) \right\}^2 \tag{9}$$

*3.2 Error back propagation learning algorithm*

Error back propagation learning algorithm has been used to update the network parameters (weights) and for minimizing error function of the ANN. For fractional differential equation we consider an unsupervised version of back propagation method. Here gradient descent method has been used for updating the parameters.

$$w_j^{k+1} = w_j^k + \Delta w_j^k = w_j^k + \left( -\eta \frac{\partial E(x,\Omega)^k}{\partial w_j^k} \right) \tag{10}$$

$$v_j^{k+1} = v_j^k + \Delta v_j^k = v_j^k + \left(-\eta \frac{\partial E(x,\Omega)^k}{\partial v_j^k}\right) \tag{11}$$

where $\eta$ is learning parameter, $k$ is iteration step which is used to update the weights as usual in ANN and $E(x,\Omega)$ is the error function.

### 3.3 Formulation of fractional order initial value problem for $\alpha \in (0,1]$

Let us consider a FDE of order $\alpha \in (0,1]$

$$(_{x_0}D_x^\alpha y(x)) = f(x, y(x)) \tag{12}$$

with initial condition $y(x_0) = y_0$

The approximate solution of ANN may be written as

$$y_N(x,\Omega) = y_0 + (x - x_0)N(x,\Omega) \tag{13}$$

The approximate solution $y_N(x,\Omega)$ satisfies the initial conditions and the error function may be computed as follows

$$E(x,\Omega) = \sum_{i=1}^{h} \{(_{x_0}D_x^\alpha y_N(x_i,\Omega)) - f(x_i, y_N(x_i,\Omega))\}^2 \tag{14}$$

### 3.4 Structure of multi-layer ANN model for FDE

We consider a three layer ANN model for the present problem. Fig. 1 depicts the structure of neural network architecture, which consists an input layer with single input node and a bias, one hidden layer having five hidden nodes and output layer consists one output node. Initial weights $w_j$ from input to hidden layer and $v_j$ from hidden to output layer are considered as random.

Architecture of the three layer ANN with five hidden nodes, single input and output layer (with one node) is shown in Fig. 1.

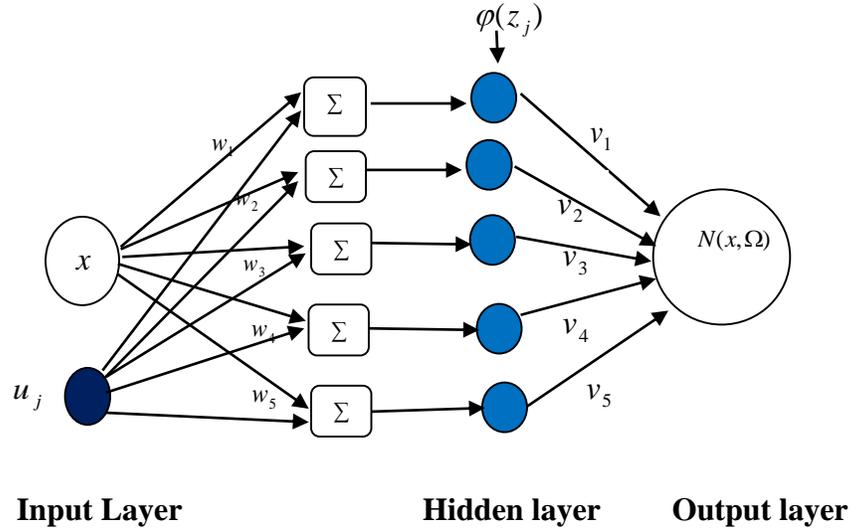

**Fig. 1** Proposed Artificial neural network architecture

*3.5 Computation of gradient for Fractional order initial value problem*

For minimizing the error function $E(x,\Omega)$ that is to update the network parameters (weights), we differentiate $E(x,\Omega)$ with respect to the parameters. Thus the gradient of network output with respect to their inputs is computed as below.

As such, the fractional derivatives (according to conformable fractional derivative) of $N(x,\Omega)$ with respect to input $x$ is written as

$$(D_x^\alpha(N(x,\Omega))) = v_j \varphi'(z_j) w_j x^{1-\alpha} \tag{15}$$

Let $(D_x^\alpha(N(x,\Omega))) = N_\beta$ denotes the derivative of the network output with respect to its inputs. The derivative of $N_\beta$ with respect to other parameters may be obtained as (according to conformable fractional derivative rules)

$$\frac{\partial N_\beta}{\partial w_j} = v_j x^{1-\alpha}(\varphi'(z_j) + \varphi''(z_j)w_j x) \tag{16}$$

$$\frac{\partial N_\beta}{\partial v_j} = w_j x^{1-\alpha}\varphi'(z_j) \tag{17}$$

$$\frac{\partial N_\beta}{\partial u_j} = w_j v_j x^{1-\alpha}\varphi''(z_j) \tag{18}$$

here

$$N(x,\Omega) = \sum_{j=1}^{m} v_j \varphi(z_j) \text{ and } z_j = w_j x + u_j$$

From (13) we have (by differentiating)

$$({}_{x_0}D_x^\alpha y_N(x)) = (x-x_0)^{1-\alpha}(N(x,\Omega)) + (x-x_0)(D_x^\alpha(N(x,\Omega))) \tag{19}$$

After simplifying the above equation we get

$$({}_{x_0}D_x^\alpha y_N(x)) = (x-x_0)^{1-\alpha}(N(x,\Omega)) + (x-x_0)(w_j v_j \varphi'(z_j)x^{1-\alpha}) \tag{20}$$

## 4. Numerical Examples and discussion

These following informative example problems are included to understand the proposed method. The approximate results by ANN model are compared with analytical\existing numerical solutions of each example to show the powerfulness of the proposed method.

**Example 1:**

Let us consider a fractional order differential equation

$$({}_0D_x^\alpha y(x)) = x \qquad x \in [0,1]$$

with initial condition $y(0) = 0$

The corresponding ANN approximate solution is expressed as

$$y_N(x, \Omega) = xN(x, \Omega)$$

The network is trained for ten equidistant points in [0, 1] and five hidden nodes. Table 1 incorporates the analytical and ANN solutions for $\alpha = 0.5$. Fig. 2 shows comparison between analytical and ANN results ($\alpha = 0.5$). Comparison between analytical and ANN solutions for $\alpha = 0.75$ and 0.85 are cited in Figs. 3 and 4 respectively. The error function for $\alpha = 0.5, 0.75$ and 0.85 are shown in Fig. 5. From Figs. 2, 3 and 4, one may see that the approximate solutions from our proposed ANN method are excellent agreement with the exact solutions for different values of $\alpha$.

**Table 1:** Analytical and ANN results for $\alpha = 0.5$ (Example 1)

| Input data | Analytical | ANN |
| --- | --- | --- |
| 0 | 0 | 0 |
| 0.1000 | 0.0238 | 0.0235 |
| 0.2000 | 0.0673 | 0.0663 |
| 0.3000 | 0.1236 | 0.1285 |
| 0.4000 | 0.1903 | 0.1900 |
| 0.5000 | 0.2660 | 0.2676 |
| 0.6000 | 0.3496 | 0.3287 |
| 0.7000 | 0.4406 | 0.4389 |
| 0.8000 | 0.5383 | 0.5481 |
| 0.9000 | 0.6423 | 0.6595 |
| 1.0000 | 0.7523 | 0.7595 |

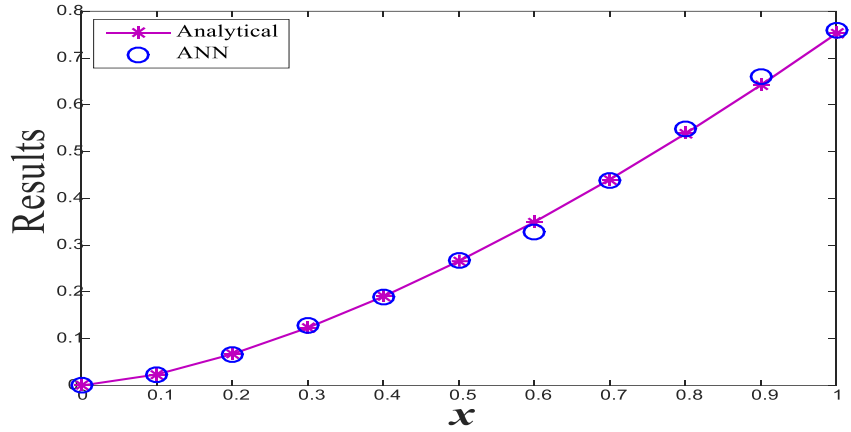

**Fig.2** Plot of Analytical and ANN results for $\alpha = 0.5$ (Example 1)

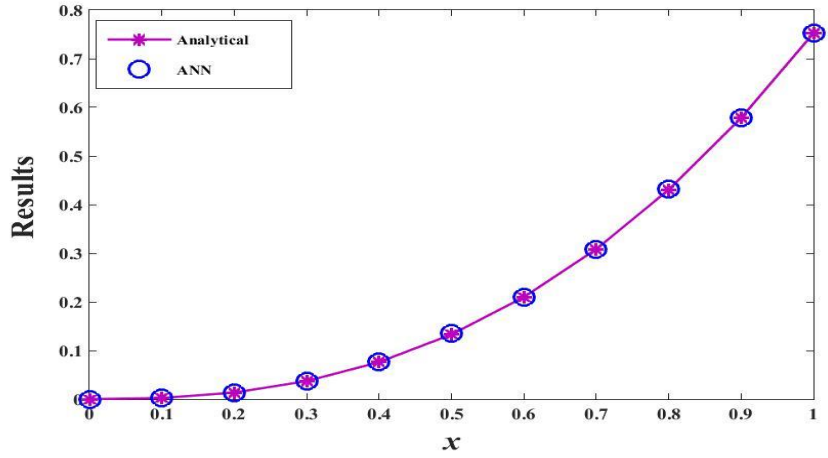

**Fig. 3** Plot of Analytical and ANN results for $\alpha = 0.75$ (Example 1)

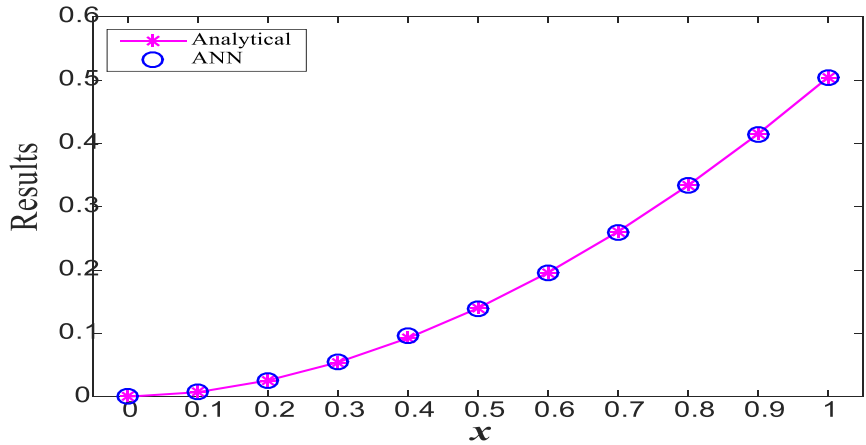

**Fig. 4** Plot of Analytical and ANN results for $\alpha = 0.85$ (Example 1)

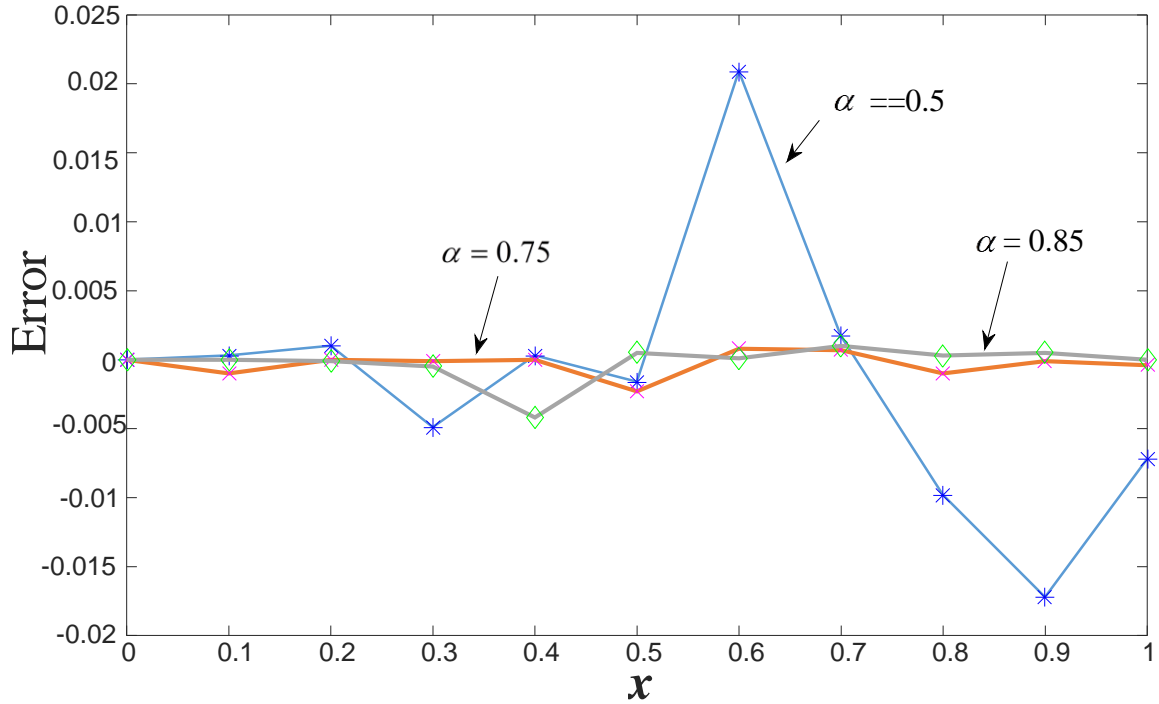

**Fig. 5** Plot of error between Analytical and ANN results for $\alpha = 0.5, 0.75$ and $0.85$ (Example 1)

**Example 2:**

Let us take now the flowing fractional order initial value problem

$$({}_0D_x^\alpha y(x)) = x^2 + 2x^{3/2} - y$$

subject to $y(0) = 0$

As discussed above we can write the ANN approximate solution as

$$y_N(x, \Omega) = xN(x, \Omega)$$

We train the network for ten equidistant points in the domain [0, 1] with five hidden nodes and $\alpha = 0.5$. Table 2 shows comparison between analytical and approximate ANN solutions for $\alpha = 0.5$. Comparison between analytical and ANN solutions are depicted in Fig.6. Error function has been plotted in Fig.7.

**Table 2:** Analytical and ANN results for $\alpha = 0.5$ (Example 2)

| Input data | Analytical | ANN |
|---|---|---|
| 0 | 0 | 0 |
| 0.1000 | 0.0100 | 0.0130 |
| 0.2000 | 0.0400 | 0.0405 |
| 0.3000 | 0.0900 | 0.0903 |
| 0.4000 | 0.1600 | 0.1608 |
| 0.5000 | 0.2500 | 0.2534 |
| 0.6000 | 0.3600 | 0.3513 |
| 0.7000 | 0.4900 | 0.4889 |
| 0.8000 | 0.6400 | 0.6469 |
| 0.9000 | 0.8100 | 0.8270 |
| 1.0000 | 1.0000 | 0.9926 |

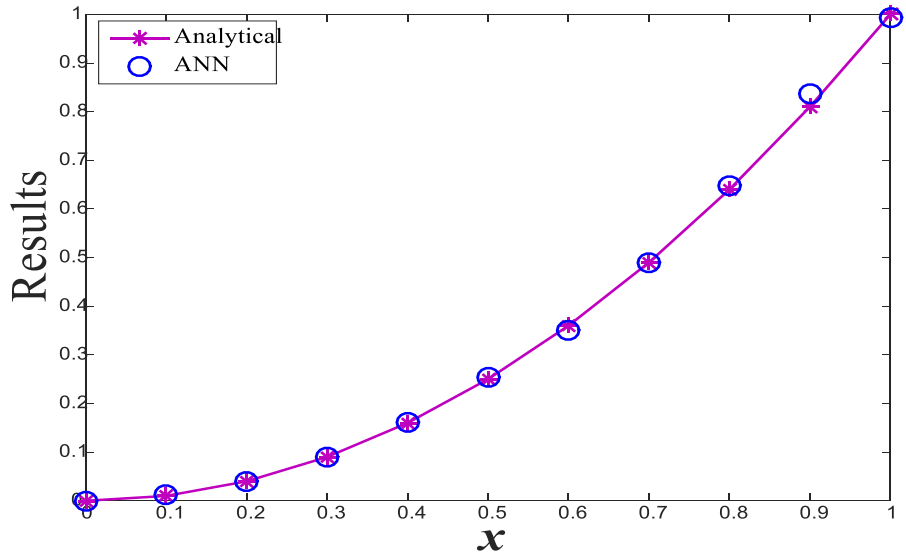

**Fig.6** Plot of Analytical and ANN results (Example 2)

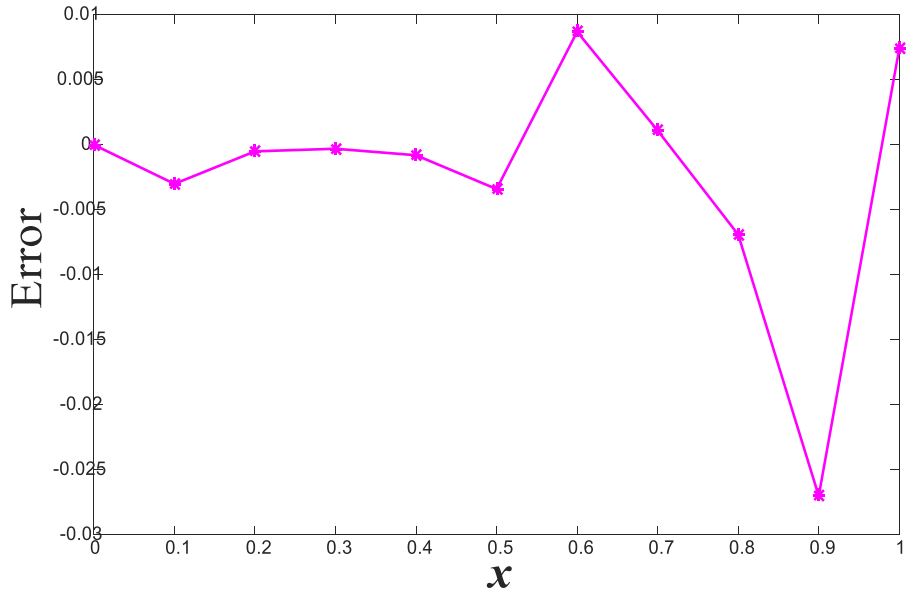

**Fig. 7** Plot of error between Analytical and ANN results (Example 2)

**Example 3:**

Finally, we take a fractional Reccati differential equation

$$(_0D_x^\alpha y(x)) = 1 - y^2 \qquad \text{for} \quad 0 < \alpha \leq 1 \text{ and } 0 \leq x < 1$$

with initial condition $y(0) = 0$

This equation has the analytical solution

$$y(x) = \frac{e^{2/\alpha\, x^\alpha} - 1}{e^{2/\alpha\, x^\alpha} + 1}$$

The ANN approximate solution in this case is represented as

$$y_N(x, \Omega) = xN(x, \Omega)$$

Again the network is trained with ten equidistant points. Table 3 incorporates the comparison among numerical solutions obtained by Chebyshev wavelet method (CWM), Ref. [37] and our proposed ANN model present for $\alpha = 0.5, 0.75$ and 1. Figs. 8 and 9 show that the approximate results from our proposed ANN method are in good agreement with the results obtained by CWM and Ref. [37] for $\alpha = 0.5, 0.75$. Comparison among the CWM, Ref. [37] and proposed ANN results for $\alpha = 1$ is illustrated in Fig.10. Plot of the error function between CWM and ANN solutions for the values of $\alpha = 0.5, 0.75$, exact and ANN for $\alpha = 1$ is cited in Fig. 11.

**Table 3:** CWM, Ref [37] and ANN results for $\alpha = 0.5, 0.75$ and 1 (Example 3)

| Input data | For $\alpha = 0.5$ | | | For $\alpha = 0.75$ | | | For $\alpha = 1$ | | |
|---|---|---|---|---|---|---|---|---|---|
| | CWM | Ref.[37] | ANN | CWM | Ref.[37] | ANN | CWM | Ref.[37] | ANN |
| 0 | 0 | 0 | 0 | 0 | 0 | 0 | 0 | 0 | 0 |
| 0.1000 | 0.3301 | 0.3214 | 0.3299 | 0.1901 | 0.1893 | 0.1911 | 0.0997 | 0.0996 | 0.0995 |
| 0.2000 | 0.4367 | 0.4326 | 0.4352 | 0.3098 | 0.3112 | 0.3102 | 0.1974 | 0.1973 | 0.1974 |
| 0.3000 | 0.5048 | 0.5178 | 0.5040 | 0.4046 | 0.4078 | 0.4116 | 0.2913 | 0.2913 | 0.2909 |
| 0.4000 | 0.5538 | 0.5540 | 0.5521 | 0.4817 | 0.4853 | 0.4837 | 0.3799 | 0.3799 | 0.3796 |
| 0.5000 | 0.5912 | 0.6088 | 0.5919 | 0.5451 | 0.5175 | 0.5390 | 0.46211 | 0.4622 | 0.4622 |
| 0.6000 | 0.6210 | 0.6381 | 0.6290 | 0.5977 | 0.5977 | 0.5929 | 0.5370 | 0.5371 | 0.5372 |
| 0.7000 | 0.6454 | 0.6464 | 0.6486 | 0.6419 | 0.6385 | 0.6401 | 0.6044 | 0.6043 | 0.6044 |
| 0.8000 | 0.6601 | 0.6584 | 0.6581 | 0.6789 | 0.6726 | 0.6609 | 0.6641 | 0.6638 | 0.6640 |
| 0.9000 | 0.6835 | 0.6589 | 0.6891 | 0.7101 | 0.7028 | 0.7121 | 0.7162 | 0.7156 | 0.7163 |

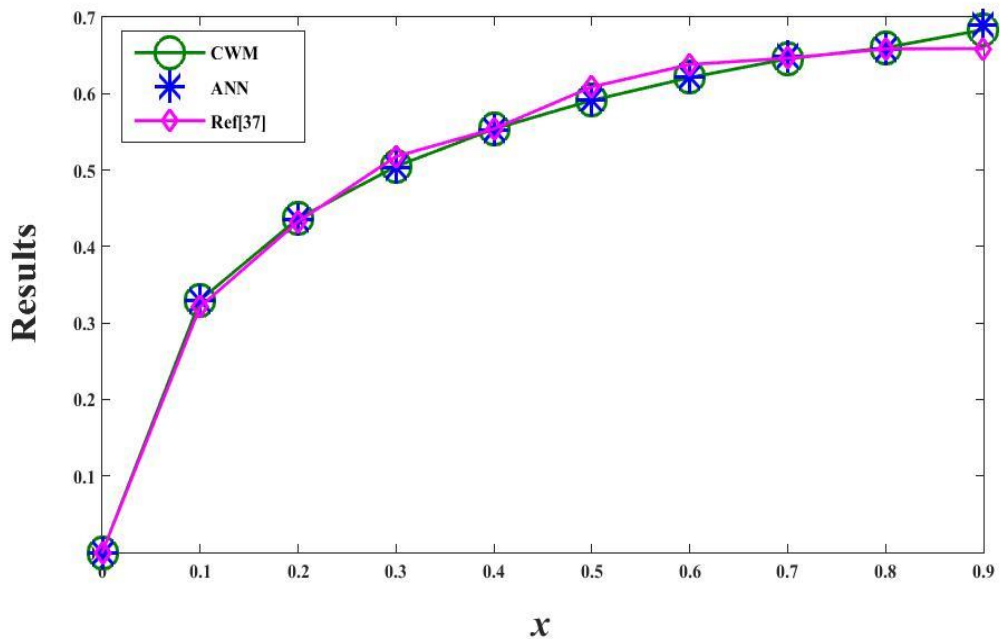

**Fig. 8** Plot of CWM, ANN and Ref. [37] results for $\alpha = 0.5$ (Example 3)

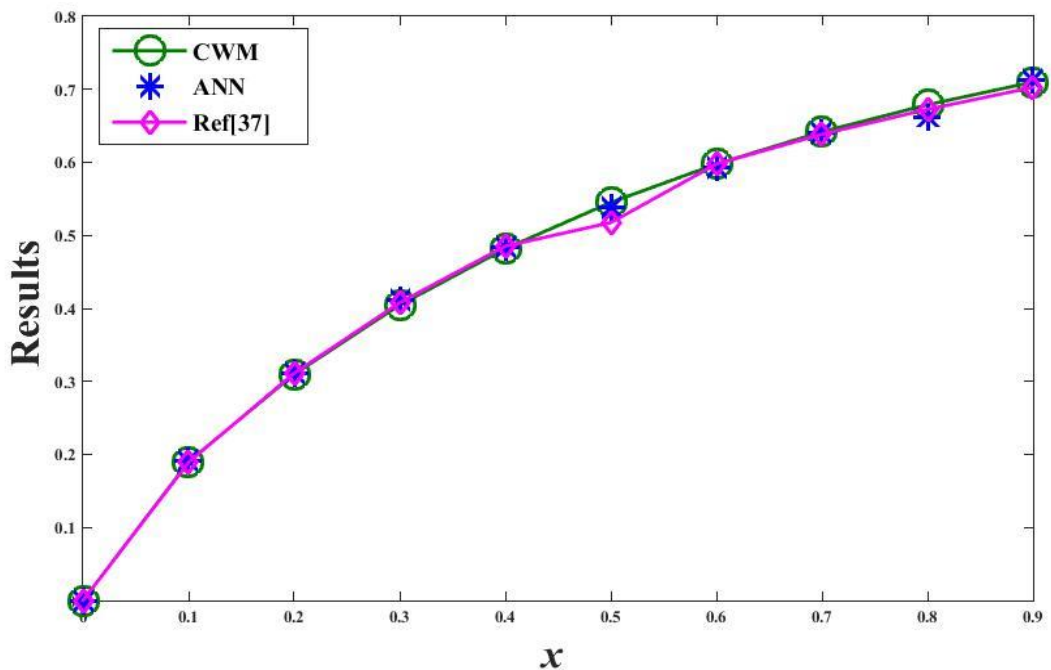

**Fig. 9** Plot of CWM, ANN and Ref. [37] results for $\alpha = 0.75$ (Example 3)

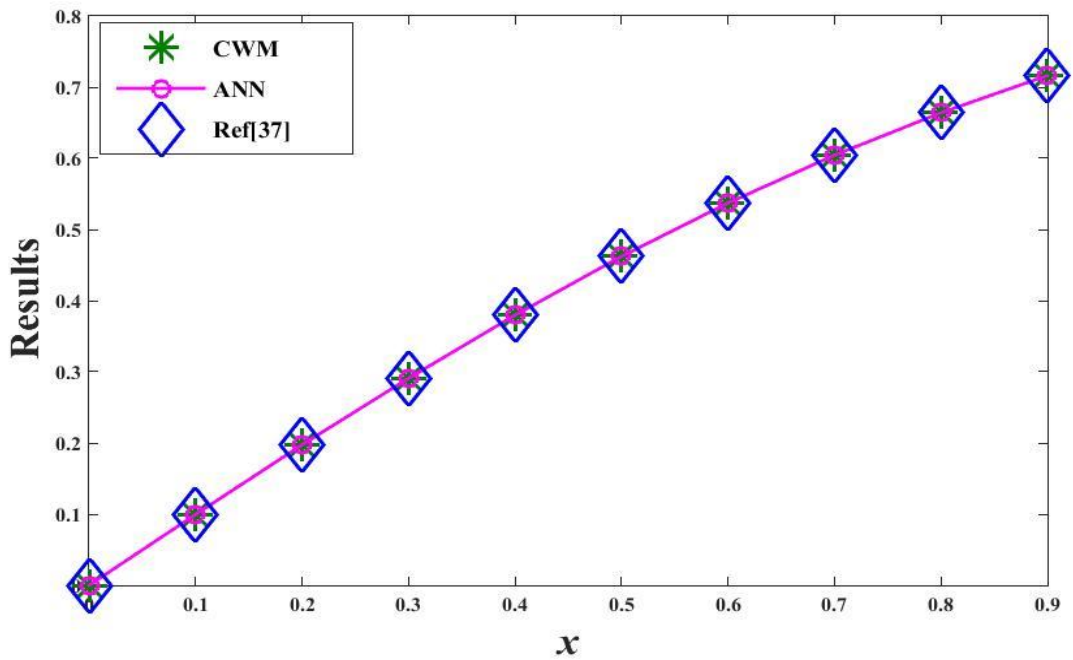

**Fig. 10** Plot of CWM, ANN and Ref. [37] results for $\alpha = 1$ (Example 3)

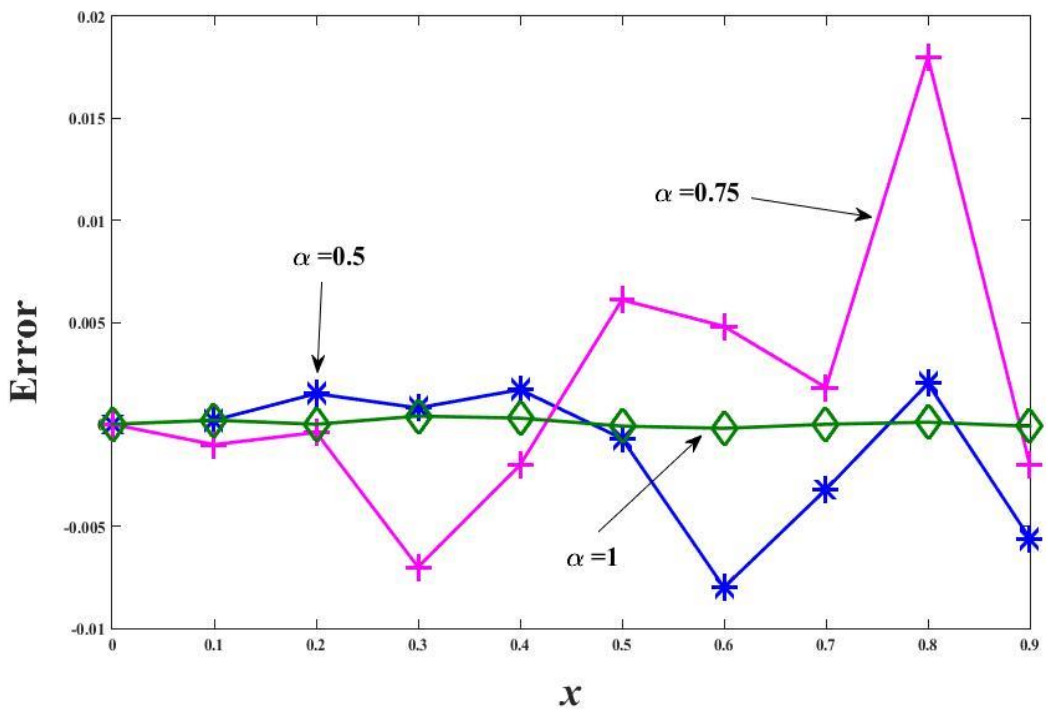

**Fig. 11** Plot of error functions for $\alpha = 0.5, 0.75$ and 1 (Example 3)

## Conclusion

This paper presents a new approach to solve fractional order differential equations by using artificial neural network model. Accuracy of the proposed method has been examined by solving various FDEs. Moreover the algorithm is unsupervised and error back propagation algorithm is used to minimize the error function. Corresponding initial weights from input to hidden and hidden to output are random. The approximate solution is closed and differentiable. One may see from the tables and graphs that the approximate solution by ANN makes the results more accurate. Lastly it may be mentioned that the ANN algorithm is simple, computationally efficient and straight forward.

## Acknowledgement

The first author is thankful to Department of Science and Technology (DST), Government of India for financial support under Women Scientist Scheme-A.